\documentclass[12pt,a4paper]{article}
\usepackage{hyperref}
\hypersetup{pdfstartview=FitH}
\hypersetup{pdfpagemode=FitWidth}
\usepackage{amssymb}
\usepackage{amsmath}
\usepackage{times} 
\newcommand{\n}{\newcommand} 
\n{\oo}{\operatorname}
\n{\mb}{\mathbb}
\n{\bb}{\bigskip}
\n{\SO}{\oo{S O}} 
\n{\Z}{\mb{Z}}
\newcommand{\e}{\equiv}
\newcommand{\f}{\varphi}
\newcommand{\ii}{\hskip 1em\relax}

\newcommand{\q}{\quad} 
\newcommand{\ttt}{\widetilde}
\newcommand{\vv}{\overline}

\begin{document}
\begin{center}
{\Huge The Gauss-Dirichlet Orbit Number}
\end{center}\bb

\centerline{\em To the memory of Pierre Kaplan}\bb

In \cite{dd} (Section 91) Dirichlet attaches a nonnegative integer $C(d,m)$ to each pair $(d,m)$ where $d$ is an integer congruent to 0 or 1 mod 4 and $m$ is a nonzero integer, and expresses, in some particular cases, $C(d,m)$ in terms of Jacobi symbols. The construction was already implicit in Article 169 of Gauss's {\em Disquisitiones Arithmeticae} (called ``{\em Disquisitiones}'' henceforth), and can be described as follows.\bb

\ii Let $F$ be the (nonempty) set of those quadratic forms 
 
$$f:=[a,b,c]:=a\,X^2+b\,X\,Y+c\,Y^2,$$

where $a,b,c$ are integers satisfying $b^2-4a c=d$, and denote again by $f$ the corresponding function from $\Z^2$ to $\Z$. Let $S$ be a representative system for the orbits of the natural action of $G:=\oo{S L}(2,\Z)$ on $F$; observe that the stabilizer $\SO(f)$ of $f\in S$ in $G$ acts on $f^{-1}(m)$; and let $C_f$ be the cardinality of the orbit set. Then $C(d,m)$ is the sum of the $C_f$ when $f$ runs over $S$:

$$C(d,m):=\sum_{f\in S}|\SO(f)\backslash f^{-1}(m)|.$$

Clearly $C(d,m)$ doesn't depend on the choice of $S$.\bb

\ii In \cite{l} (Theorem 204) Landau enlarges the validity domain of Dirichlet's formula by replacing Jacobi symbols with Kronecker symbols. The integer $C(d,m)$ is computed for $d$ nonsquare and $m$ arbitrary by Huard, Kaplan and Williams in \cite{hkw} (Theorem 9.1), and by Sun and Williams in \cite{sw} (Theorems 4.1 and 4.2). In \cite{sw} the authors make the crucial observation that $C(d,m)$ is multiplicative in $m$. They prove this multiplicativity for nonsquare $d$, but it holds for all $d$, as follows immediately from Article 169 of the {\em Disquisitiones}. In \cite{r} Rudnick gives a definition of $C(d,m)$ which is clearer than Dirichlet's. Here we compute $C(d,m)$ in full generality. Thank you to Keith Matthews, Ze\'ev Rudnick and Kenneth Williams for their friendly and efficient help.\bb

{\bf Gauss's Observation.} For $f$ in $S$, let $f^{-1}(m)'$ be the set of those elements of $f^{-1}(m)$ which have coprime coordinates; and put 

$$T:=\{t\in\Z/2m\Z\ |\ t^2\equiv d\bmod 4m \}.$$ 

\ii Then there is a unique map $\f_f$ from $f^{-1}(m)'$ to $T$ such that for any $z$ in $f^{-1}(m)'$ and any $g$ in $G$ mapping $(1,0)$ to $z$ we have 

$$f g=[m,n,\ell],\q\f_f(z)=n\bmod2m.$$ 

\ii Moreover the $\f_f$ induce a bijection 

$$\coprod_{f\in S}\SO(f)\backslash f^{-1}(m)'
\stackrel{\sim}{\longrightarrow}T.$$

In particular $C(d,m)$ is finite. Moreover we have $C(d,m n)=C(d,m)C(d,n)$ whenever $m$ and $n$ are coprime --- property which we express by saying that the function $C(d,?)$ is {\bf multiplicative}. \bb 

\ii Here is a proof of Gauss's Observation. Put 
 
$$\ttt{T}:=\{t\in\Z\ |\ t^2\equiv d\bmod 4m \}.$$ 

\ii Define, for $f$ in $S$, the map $\ttt{\f}_f$ from $f^{-1}(m)'$ to $\ttt{T}$ as follows. Let $z$ be in $f^{-1}(m)'$. There is a $g$ in $G$ which maps $(1,0)$ to $z$. Choose such a $g$, note that $f g$ is equal to $[m,n,\ell]$ for some $(n,\ell)$ in $\ttt{T}\times\Z$, and define $\ttt{\f}_f(z)$ as $n$.\bb

\ii Define the map $\ttt{\psi}$ from $\ttt{T}$ to $\coprod f^{-1}(m)'$ as follows. Let $n$ be in $\ttt{T}$. Then there is a unique pair $(\ell,f)$ in $\Z\times S$ for which $[m,n,\ell]$ is equal to $f g$ for some $g$ in $G$. Choose such a $g$, note that $g(1,0)$ is in $f^{-1}(m)'$, and define $\ttt{\psi}(n)$ as $g(1,0)$.\bb

\ii One checks that the map 

$$\f:\coprod f^{-1}(m)'\to T$$ 

induced by the $\ttt{\f}_f$ doesn't depend on the choices made to define the $\ttt{\f}_f$, that the map 

$$\psi:\ttt{T}\to \coprod\SO(f)\backslash f^{-1}(m)'$$ 

induced by $\ttt{\psi}$ doesn't depend on the choices made to define $\ttt{\psi}$, that $\f$ factors through a map 

$$\vv{\f}:\coprod\SO(f)\backslash f^{-1}(m)'\to T,$$ 

that $\psi$ factors through a map 

$$\vv{\psi}:T\to\coprod\SO(f)\backslash f^{-1}(m)',$$ 

and that $\vv{\f}$ and $\vv{\psi}$ are inverse.\bb

\ii Of course there are many verifications to make, but they are straightforward. Obviously Gauss couldn't use exactly this language, but the substance of the above argument is clearly contained in Article 169 of the {\em Disquisitiones}.\bb

\ii Let us fix $d$ and let $m$ vary. To prove the multiplicativity of $C(d,m)$ it suffices to prove that of the cardinality of $T$. We have the following lemma. Let $\alpha$ and $\mu$ be arithmetic functions defined on the nonzero integers. Assume $\mu$ is multiplicative, $2\alpha(x)=\mu(4x)$ for all $x$, and $\mu(4)=2$. Then $\alpha$ is multiplicative. Now take the cardinality of $T$ as $\alpha(m)$, and take the number of square roots of $d$ mod $x$ as $\mu(x)$ --- a multiplicative function by the Chinese Remainder Theorem.\bb

{\bf Definition of $K(d,m)$.} Fix an integer $d$ congruent to 0 or 1 mod 4. Let $\chi$ be the map from the positive integers to the integers characterized by the following properties: \bb 

$\bullet$ $\chi(1)=1$, $\chi(p)=0$ if $p$ is a prime dividing $d$, \smallskip

$\bullet$ $\chi(p)$ is the Legendre symbol $(\frac{d}{p})$ if $p$ is an odd prime not dividing $d$,\smallskip

$\bullet$ $\chi(2)=1$ if $d\equiv1$ mod 8, $\chi(2)=-1$ if $d\equiv5$ mod 8,\smallskip

$\bullet$ $\chi(m n)=\chi(m)\chi(n)$ for all $m,n$. \bb 

For $n>0$ define the {\bf Kronecker symbol} by $(\frac{d}{n}):=\chi(n)$, and for all nonzero integer $m$ put 

$$K(d,m):=\sum_{0<\mu|m}\left(\frac{d}{\mu}\right).$$

In particular $K(d,?)$ is multiplicative.\bb

{\bf Claims.} If $d$ and $m$ are coprime then $C(d,m)=K(d,m)$. We clearly have $C(d,\pm1)=1$. By multiplicativity it only remains to compute $C(d,m)$ when $m$ is a power of a prime factor of $d$. Fix such a prime factor $p$, and write $B(d,n)$ for $C(d,p^n)$, where $n$ is a nonnegative integer. As $B(d,0)$ and $B(d,1)$ are equal to 1 (easy verification), we henceforth assume $n\ge2$.\bb

\ii Putting $k:=\left\lfloor n/2\right\rfloor+1$ we get $B(0,n)=(p^k-1)/(p-1)$.\bb 

\ii Let $p$ be odd and $d$ a nonzero multiple of $p$. Write $d$ as $p^\delta D$ with $\delta>0$ and $D$ prime to $p$. As we clearly have $B(d,n)=B(0,n)$ --- which is already computed --- for $n\le\delta$, we can assume $n>\delta$. Recall that $d$ is the square of a $p$-adic integer if and only if $\delta$ is even and the Legendre symbol $(\frac{D}{p})$ equals $1$. Put $k:=\lceil (n-\delta)/2\rceil$. Then $B(d,n)$ is equal to $B(0,n-2k)+2k p^{\delta/2}$ if $d$ is the square of a $p$-adic integer, and to $B(0,n-2k)$ otherwise.\bb

\ii Consider now the case $p=2$. Let $\delta$ be a positive integer.\bb

\ii Assume $d=2^{2\delta+1} D$ with $D$ odd. Then $B(d,n)$ is equal to $B(0,n)$ if $n\le 2\delta-1$, and to $2^\delta-1$ if $n\ge2 \delta$.\bb

\ii Assume $d=2^{2\delta} D$ with $D\e3$ mod 4. Then $B(d,n)$ is equal to $B(0,n)$ if $n\le2\delta-2$, and to $2^\delta-1$ if $n\ge2\delta-2$.\bb 

\ii Assume $d=2^{2\delta} D$ with $D\e5$ mod 8. Then $B(d,n)$ is equal to $B(0,n)$ if $n\le2\delta-2$. Suppose $n\ge2\delta-1$. Then $B(d,n)$ is equal to $2^\delta-1$ if $n$ is odd, and to  $2^{\delta+1}-1$ if $n$ is even.\bb

\ii Assume $d=2^{2\delta}D$ with $D\e1$ mod 8, {\em i.e.} $d$ is the square of a 2-adic integer. Then $B(d,n)$ is equal to $B(0,n)$ if $n\le2\delta-2$. For $k\ge0$ we have 

$$B(d,2\delta-1+2k)=2^\delta(2k+1)-1,\q 
B(d,2\delta+2k)=2^{\delta+1}(k+1)-1.$$\smallskip

{\bf Corollary 1.} The functions $C(d,?)$ and $K(d,?)$ coincide if and only if $d$ has no odd square factor and $d$ is not an even square mod 16.\bb 

{\bf Corollary 2.} The function $n\mapsto C(d,p^n)$ is bounded if and only if $d$ is not the square of a $p$-adic integer.\bb 

{\bf Lemma.} For $d,n$ as above, let $2B'(d,n)$ be the number of roots of the congruence 

$$x^2\e d\bmod4p^n.$$

In particular, if $p>2$, then $B'(d,n)$ is the number of roots of $x^2\e d\bmod p^n$. Moreover we have

$$B(d,n)=\sum_{0\le j\le \lfloor n/2\rfloor}B'(d,n-2j).$$

{\bf Proofs.} Gauss's Observation, along with Articles 104 and 105 of the {\em Disquisitiones}, implies the Lemma and the Claims.\bb

\ii In the following two examples $S$ reduces to a singleton. See the {\em Disquisitiones} or \cite{dd} for a proof of this fact. The other statements follow from Corollary~1. Let $m$ be a nonzero integer.\bb 

{\bf Example 1}. Let $d$ be $-4$. We can take $S:=\{f\}$ with $f:=X^2+Y^2$. We have $|\SO(f)|=4$, and the number of representations of $m$ as a sum of two squares is four times $M-N$, where $M$ is the number of divisors of $m$ congruent to 1 mod 4, whereas $N$ is the number of divisors of $m$ congruent to $-1$ mod 4. \bb 

{\bf Example 2}. Let $d$ be 8. We can take $S:=\{f\}$ with $f:=X^2-2Y^2$. Then $C(8,m)$ is the number of divisors of $m$ congruent to $\pm1$ mod 8 minus the number of divisors of $m$ congruent to $\pm3$ mod 8.

\vfill Pierre-Yves Gaillard\hfill\tiny gauss-dirichlet.orbit.number.080910, Wed Sep 10 08:55:21 CEST 2008. 

\end{document}